\documentclass[11pt]{article}%
\usepackage[latin1]{inputenc}
\usepackage{latexsym}
\usepackage{amsthm}
\usepackage[Sonny]{fncychap}
\usepackage{enumerate}
\usepackage[dvips]{color,graphicx}
\usepackage{pst-all}
\usepackage{pstcol}
\usepackage{amsbsy}
\usepackage{amsmath}
\usepackage{amsfonts}
\usepackage{amssymb}
\usepackage{graphicx}
\usepackage{amscd}%
\newtheorem{definicion}{Definition}[section]

\newtheorem{theorem}[definicion]{Theorem}
\newtheorem{corollary}[definicion]{Corollary}
\newtheorem{lemma}[definicion]{Lemma}
\newtheorem{remark}[definicion]{Remark}

\numberwithin{equation}{section}

\def\qed { \vskip 0pt \hfill \hbox{\vrule height 5pt width 5pt depth 0pt} \vskip 10pt}
\begin{document}

\title{An approach without using Hardy inequality for the linear heat equation with singular potential}
\author{{{Lucas C. F. Ferreira} {\thanks{L.C.F. Ferreira was supported by FAPESP and
CNPQ, Brazil. (corresponding author)}}}\\{\small Universidade Estadual de Campinas, IMECC - Departamento de
Matem\'{a}tica,} \\{\small {Rua S\'{e}rgio Buarque de Holanda, 651, CEP 13083-859, Campinas-SP,
Brazil.}}\\{\small \texttt{Email:lcff@ime.unicamp.br}} \vspace{0.5cm}\\{{Cl\'{a}udia Aline A. S. Mesquita}{\thanks{C.A.A.S. Mesquita was supported by
CNPQ, Brazil.}}}\\{\small Universidade Estadual de Campinas, IMECC - Departamento de
Matem\'{a}tica,} \\{\small {Rua S\'{e}rgio Buarque de Holanda, 651, CEP 13083-859, Campinas-SP,
Brazil.}}\\{\small \texttt{Email:\ ra098154@ime.unicamp.br}}}
\date{}
\maketitle

\noindent\textsc{abstract}. The aim of this paper is to employ a strategy
known from fluid dynamics in order to provide results for the linear heat
equation $u_{t}-\Delta u-V(x)u=0$ in $\mathbb{R}^{n}$ with singular
potentials. We show well-posedness of solutions, without using Hardy
inequality, in a framework based in the Fourier transform, namely $PM^{k}%
$-spaces. For arbitrary data $u_{0}\in PM^{k}$, the approach allows to compute
an explicit smallness condition on $V$ for global existence in the case of $V$
with finitely many inverse square singularities. As a consequence,
well-posedness of solutions is obtained for the case of the monopolar
potential $V(x)=\frac{\lambda}{\left\vert x\right\vert ^{2}}$ with $\left\vert
\lambda\right\vert <\lambda_{\ast}=\frac{(n-2)^{2}}{4}$. This threshold value
is the same one obtained for the global well-posedness of $L^{2}$-solutions by
means of Hardy inequalities and energy estimates. Since there is no any
inclusion relation between $L^{2}$ and $PM^{k}$, our results indicate that
$\lambda_{\ast}$ is intrinsic of the PDE and independent of a particular
approach. We also analyze the long-time behavior of solutions and show there
are infinitely many possible asymptotics characterized by the cells of a
disjoint partition of the initial data class $PM^{k}$.

\


\


\newpage

\section{Introduction and statements of results}

\

\bigskip We concern with the Cauchy problem for the liner heat equation
\begin{align}
u_{t}-\Delta u-V(x)u  &  =0\text{\hspace{1.2cm} in \ }\mathbb{R}%
^{n}\label{heat1}\\
u(x,0)  &  =u_{0}(x)\text{\hspace{0.5cm} in \ }\mathbb{R}^{n}, \label{heat3}%
\end{align}
where $n\geq3$ and $V(x)$ is a singular potential. Of particular interest are
the negative power law ones which appear in a number of physical and
mathematical contexts, such as molecular physics, non-relativistic quantum
mechanics, quantum cosmology, linearized analysis of combustion models, and
many others (see \cite{FMT1},\cite{FMT2},\cite{Frank},\cite{Landau}%
,\cite{Levy},\cite{Peral-Vazquez} and references therein). These potentials
can be classified according to the number of singularities (poles), $\sigma
$-degree of the singularity (order of the poles), dependence on directions
(anisotropy) and decay at infinity.

A critical situation is when $\sigma$ coincides with the order of the main
part of the associated elliptic operator, e.g. $\sigma=2$ for (\ref{heat1}%
)-(\ref{heat3}). This case presents further difficulties in its mathematical
analysis due to the following features: $V$ does not belong to Kato's class
and $Vu$ cannot be handled as a lower order term. Examples of those are the
inverse square (Hardy) potential
\begin{equation}
V(x)=\frac{\lambda}{\left\vert x\right\vert ^{2}},\text{ for }\lambda
\in\mathbb{R},\label{pot1}%
\end{equation}
and
\begin{equation}
V(x)=\sum_{j=1}^{m}\frac{\lambda_{j}}{|x-x^{j}|^{2}}\text{ and }%
V(x)=\sum_{j=1}^{m}\frac{(x-x^{j}).d^{j}}{|x-x^{j}|^{3}}\text{,}\label{pot2}%
\end{equation}
\ where $x^{j}=(x_{1}^{j},x_{2}^{j},...,x_{n}^{j})\in$\ $\mathbb{R}^{n}$ and
$d^{j}\in\mathbb{R}^{n}$ are given constant vectors. The potentials in
(\ref{pot2}) are called isotropic and anisotropic multipolar inverse square
ones, respectively.

There is a well-known concept of criticality associated to the size of the
parameter $\lambda$ in (\ref{pot1}) with respect to best constant in Hardy's
inequality $\lambda_{\ast}=\frac{(n-2)^{2}}{4}$, which reads as
\begin{equation}
\lambda_{\ast}\int_{\mathbb{R}^{n}}\frac{u^{2}}{\left\vert x\right\vert ^{2}%
}dx\leq\left\Vert \nabla u\right\Vert _{L^{2}(\mathbb{R}^{n})}^{2}.
\label{Hardy1}%
\end{equation}
In fact, the celebrated work of Baras and Goldstein \cite{Baras1} established
a threshold value that decides whether (or not) positive solutions in
$L^{2}(\mathbb{R}^{n})$ exist. Precisely, if $0\leq\lambda\leq\lambda_{\ast}$
then (\ref{heat1})-(\ref{heat3}) is well-posedness in $L^{2}(\mathbb{R}^{n}),$
and it is not well-posedness for $\lambda>\lambda_{\ast}$. In the latter case,
there is no nontrivial nonnegative solutions for $u_{0}\geq0$, while weak
positive solutions do exist in the former one when $u_{0}\geq0$ and
$u_{0}\not \equiv 0.$ Because of this dichotomy, the cases $\lambda\in
\lbrack0,\lambda_{\ast}),$ $\lambda=\lambda_{\ast}$ and $\lambda>\lambda
_{\ast}$ are named respectively as sub-critical, critical and supercritical
values for $\lambda$ (see \cite{Goldstein2}, \cite{Vazquez1} for a deeper discussion).

Over the last fifteen years, the results of \cite{Baras1} have motivated many
works concerning heat equations with singular potentials. As well as
\cite{Baras1}, Hardy inequality and its versions play a crucial role in the
results of the literature. In what follows, without making a complete list, we
review some important works. In a smooth bounded domain $\Omega$ and for
general positive singular potential $V\in L_{loc}^{1}(\Omega)$, existence and
non-existence results were proved in \cite{Cabre-Martel1} via conditions on
the infimum of the spectrum of the operator $\Delta-V$. Kombe \cite{Kombe}
showed that the nonexistence result of positive solutions in
\cite{Cabre-Martel1} is not affected\ when the potential $V$ is perturbed by a
highly oscillating singular sign-changing potential. Another extension for a
parabolic equation with variable coefficients in the principal part can be
found in \cite{Goldstein2}. The authors of \cite{Vazquez1} proved an
improvement of the Hardy-Poincar\'{e} inequality in bounded domains and a
weighted version of that inequality in $\mathbb{R}^{n}$. Afterwards they
showed exponential stabilization towards a solution in separated variables in
bounded domains, and polynomial stabilization towards a radially symmetric
solution in self-similar variables. Using Hardy type inequalities of
\cite{Vazquez1}, comparison results were obtained in \cite{Davila} for linear
elliptic and parabolic equations with $V\in L_{loc}^{1}(\Omega)$ positive.
Also, using some improved forms of Hardy-Poincar\'{e} inequalities and
Carleman estimates, inverse source problems have been considered for
(\ref{heat1}) with $0\leq\lambda\leq\lambda_{\ast}$ (see \cite{Vancostenoble1}
and their references). Motivated by instantaneous blowing-up of non-negative
$H_{0}^{1}(\Omega)$-solutions when $\lambda=\lambda_{\ast},$ the authors of
\cite{Chaudhuri} characterized some kinds of perturbations of the critical
potential $\frac{\lambda_{\ast}}{\left\vert x\right\vert ^{2}}$ for obtaining
existence and non-existence of $H^{1}$-solutions (see \cite{Gkikas1} for
related results in the stationary case). Nonexistence for $\lambda
>\lambda_{\ast}$ for a perturbation of (\ref{heat1}) with $-\nabla\sigma
\cdot\nabla u$ $+\frac{\lambda}{\left\vert x\right\vert ^{2}}u$ in place of
$Vu$ was addressed in \cite{Goldstein1}. Results concerning existence,
non-existence, Fujita exponent, self-similarity, bifurcations, instantaneous
blow-up for perturbations of (\ref{heat1}) by $u^{p}$ and $\left\vert \nabla
u\right\vert ^{p}$ can be found in \cite{Abdellaoui-Peral2},
\cite{Abdellaoui-Peral3}, \cite{Abdellaoui-Peral4}, \cite{Chaves},
\cite{Karachalios}, \cite{Liskevich1} , \cite{Reyes1} (see also \cite{Peral1}%
). Related linear and semilinear elliptic problems \cite{Abdellaoui-Peral1}%
,\cite{Chaves},\cite{Felli1} ,\cite{FMT1},\cite{FMT2},\cite{Smets} have been
considered with results also presenting a dichotomy due to influence of
Hardy potentials. See \cite{Fer-Mesq} for existence results for some
semilinear elliptic equations with multipolar potentials (\ref{pot2}) without
employing Hardy inequalities.

From the above works, the use of Hardy inequality (\ref{Hardy1}) imposes the
$L^{2}$-framework for $u$ and the condition $0\leq\lambda\leq\lambda_{\ast}$
for well-posedness of $L^{2}$-solutions. So, a natural question arises: does
there exist a framework different from $L^{2}$ in which (\ref{heat1}%
)-(\ref{heat3}) with (\ref{pot1}) is well-posedness for $0\leq\lambda
\leq\lambda_{\ast}$ or at least for $0\leq\lambda<\lambda_{\ast}$? In this
paper we give a positive response for this question by using $PM^{k}$-spaces
and a strategy based on Fourier transform which does not use Hardy inequality
(see Theorem \ref{t.i}). Since there is no any inclusion relation between
$L^{2}$ and $PM^{k}$, and the corresponding techniques are of different
natures, our results indicate that $\lambda_{\ast}$ is intrinsic of the PDE
(\ref{heat1}) and independent of a particular approach. Another goal is to
obtain results for multipolar potentials like (\ref{pot2}) providing explicit
conditions on parameters of the potential for global well-posedness of
solutions (see Corollary \ref{c.i}). Let us observe that, as well as in
\cite{Vazquez1}, we also consider sign-changing initial data, solutions and potentials.

From another viewpoint, parabolic regularization does not work for nontrivial
solutions when $V$ is singular, and then solutions are not smooth in
$\mathbb{R}^{n}$. This also motivates to analyze the well-posedness of
(\ref{heat1}) in a singular framework like $PM^{k}$ which contains functions
that can be strongly rough and not to decay as $\left\vert x\right\vert
\rightarrow\infty$. Indeed, Theorem \ref{t.i} and Corollary \ref{c.i} below
show that the formal semigroup associated with (\ref{heat1}) can be extended
to the $PM^{k}$-framework for a wide class of singular potential, namely $V\in
PM^{n-2}$ satisfying (\ref{pot-cond}).

The authors of \cite{Vazquez1} showed that certain good solutions with
positive data in the weighted space $L^{2}(\mathbb{R}^{n},\exp(\left\vert
x\right\vert ^{2}/4)dx)$ converges in the norm $t^{1/2}\left\Vert
\cdot\right\Vert _{L^{2}}$ for an explicitly given non-stationary solution (up
to a constant) as $t\rightarrow\infty$. Here we obtain a partition of $PM^{k}$
into infinite pairwise disjoint subsets (induced by an equivalent relation)
such that the long-time behavior of $u$ depends on the subset that contains
the data $u_{0}$ (see Theorem \ref{t.iii} and Remark \ref{Rem-Asymp}). In
particular, for $u_{0}$ belonging to certain special subsets, the solution
converges towards a stationary state given explicitly. These results show that
the asymptotic behavior of solutions in the $PM^{k}$-framework is more complex
than that for $L^{2}$-solutions, at least in an approximation $o(1)$,
as $t\rightarrow\infty$.

The problem is formulated via a functional equation obtained by formally
applying the Fourier transform in (\ref{heat1})-(\ref{heat3}) and then using
Duhamel principle. This approach has already been used in the context of fluid
mechanics and semilinear parabolic equations without singular potentials (see
e.g. \cite{Biler-Karch}, \cite{Cannone1}, \cite{Car-Fer1}, \cite{Le-Jan},
\cite{Miao1}). Nonetheless, to our knowledge, the present paper seems to be
its first application on a global existence problem for a parabolic PDE with
an optimal threshold value explicitly characterized by another technique.

The authors of \cite{Galaktionov1} studied the higher-order parabolic
equation
\begin{equation}
u_{t}+(-\Delta)^{m}u+\frac{\lambda}{\left\vert x\right\vert ^{\sigma}%
}u=0\label{heat-m}%
\end{equation}
with critical singularity $\sigma=2m$ and $n>2m.$ Among others, they extended
Baras-Goldstein results by obtaining an explicit threshold value $\lambda_{m}$
for well-posedness of $L^{2}$-solutions, which is the so-called Hardy's best
constant of multiplicative inequalities involving $V(x)=\lambda\left\vert
x\right\vert ^{-2m}.$ With a slight adaptation on the proofs, the results of
the present paper could be extended for (\ref{heat-m}) by obtaining
well-posedness of $PM^{k}$-solutions for $0\leq\left\vert \lambda\right\vert
<\lambda_{m}.$

In what follows we describe precisely our results. For each $k\geq0,$ the
space $PM^{k}$ is defined by
\begin{equation}
PM^{k}\equiv\{u\in\mathcal{S^{\prime}}:\hat{u}\in L_{loc}^{1}(\mathbb{R}%
^{n}),\text{ ess}\sup_{\xi\in\mathbb{R}^{n}}|\xi|^{k}|\hat{u}(\xi)|<\infty\},
\label{space1}%
\end{equation}
which is a Banach space with the norm
\[
\left\Vert u\right\Vert _{PM^{k}}\equiv ess\sup_{\xi\in\mathbb{R}^{n}}%
|\xi|^{k}|\hat{u}(\xi)|<\infty.
\]
Here $\widehat{\cdot}$ stands for the Fourier transform in $\mathcal{S^{\prime
}}(\mathbb{R}^{n})$ which is an extension of the classical definition $\hat
{f}(\xi)=\int_{\mathbb{R}^{n}}f(x)e^{-2\pi i\xi\cdot x}dx$ in $\mathcal{S}%
(\mathbb{R}^{n})$.\bigskip

The problem (\ref{heat1})-(\ref{heat3}) is formally equivalent to functional equation%

\begin{equation}
u(t)=G(t)u_{0}+L_{V}(u)(t),\text{ } \label{mild}%
\end{equation}
where the operators $G(t)u_{0}$ and $L_{V}(u)(t)$ are defined via Fourier
transform as%

\begin{align}
\widehat{G(t)u_{0}}(\xi,t)  &  =e^{-4\pi^{2}\left\vert \xi\right\vert ^{2}%
t}\hat{u}_{0}\label{Group}\\
\widehat{L_{V}(u)}(\xi,t)  &  =\int_{0}^{t}e^{-4\pi^{2}\left\vert
\xi\right\vert ^{2}(t-s)}\left(  \hat{V}\ast\hat{u}\right)  (\xi,s)\,ds.
\label{Tv}%
\end{align}
Notice that $G(t)$ is a convolution operator with Gaussian kernel
$g(x,t)=(4\pi t)^{-\frac{n}{2}}e^{-|x|^{2}/4t}$ (heat semigroup) and if $u$
and $V$ are regular enough, then $L_{V}(u)(t)=\int_{0}^{t}G(t-s)(Vu)(s)\,ds$
in the Bochner sense in $PM^{k}.$ However, for general $u$ and $V,$ the
operator $L_{V}(u)$ cannot be understood in such a sense, and the integral
with respect to $s$ in (\ref{Tv}) should be meant in a weak-sense like e.g.
\cite{Yamazaki} (see \cite{Biler-Karch}, \cite{Cannone1}).

Global-in-time solutions $u(x,t)$ will be sought in the time-dependent Banach
space%
\begin{equation}
X_{k}=BC_{w}\left(  [0,\infty);PM^{k}\right)  \text{ } \label{space2}%
\end{equation}
with norm $\left\Vert u\right\Vert _{X_{k}}=\sup_{t>0}\left\Vert
u(\cdot,t)\right\Vert _{PM^{k}}.$ Here $BC_{w}$ stands for the set of bounded
functions from a interval into a Banach space that are time-weakly continuous
in the sense of $\mathcal{S^{\prime}}(\mathbb{R}^{n})$ at each $t\geq0$.

Our well-posedness result reads as follows.

\begin{theorem}
\label{t.i}\noindent Suppose that $V\in PM^{n-2}$ and $u_{0}\in PM^{k}$ with
$2<k<n.$ Let $K(\theta_{1},\theta_{2},n)=(\nu_{\theta_{1}}\nu_{\theta_{2}}%
\nu_{n-\theta_{1}-\theta_{2}})/(\nu_{\theta_{1}+\theta_{2}}\nu_{n-\theta_{1}%
}\nu_{n-\theta_{2}}),$ where $\nu_{\theta}=\pi^{-\theta/2}\Gamma(\theta/2)$
and $\Gamma$ is the Gamma function.

\begin{itemize}
\item[(i)] (Existence and uniqueness) Let $C_{n-2,k}=\frac{K(2,n-k,n)}%
{4\pi^{2}}$ and assume that
\begin{equation}
\Vert V\Vert_{PM^{n-2}}<\frac{1}{C_{n-2,k}}. \label{pot-cond}%
\end{equation}
Then the functional equation (\ref{mild}) has a unique solution $u$ in
$X_{k}.$

\item[(ii)] (Hardy potential) For $V(x)=\frac{\lambda}{|x|^{2}},$ the
condition on $V$ becomes equivalent to $|\lambda|<(k-2)(n-k).$ Notice that the
maximum of $(k-2)(n-k)$ is $\lambda_{\ast}=\frac{(n-2)^{2}}{4}$ which is
reached at $k=\frac{n+2}{2}$. Then the item (i) provides a global solution $u$
for (\ref{mild}) for all $u_{0}\in PM^{1+\frac{n}{2}}$ and $0\leq\left\vert
\lambda\right\vert <\lambda_{\ast}$.

\item[(iii)] (Continuous dependence) The data-solution map $(u_{0}%
,V)\rightarrow u$ is Lipschitz continuous from $PM^{k}\times PM^{n-2}$ to
$X_{k}.$ More precisely, if $u$ and $v$ are solutions obtained in item (i)
corresponding to $u_{0},V$ and $\ v_{0},W$, respectively, then
\[\hspace{-0.7cm}
\left\Vert u-v\right\Vert _{X_{k}}\leq\frac{1}{1-C_{n-2,k}\left\Vert
V\right\Vert _{PM^{n-2}}}\left(  \Vert u_{0}-v_{0}\Vert_{PM^{k}}%
+\frac{C_{n-2,k}\Vert v_{0}\Vert_{PM^{k}}}{1-C_{n-2,k}\left\Vert W\right\Vert
_{PM^{n-2}}}\left\Vert V-W\right\Vert _{PM^{n-2}}\right)
\]

\end{itemize}
\end{theorem}

Theorem \ref{t.i} can be applied for other types of singular potentials giving
explicit conditions on the size of them.

\begin{corollary}
\label{c.i} Under hypotheses of Theorem \ref{t.i}. Let $\beta(\cdot,\cdot)$
stands for the Beta function and $\left\vert \cdot\right\vert $ the sum norm
in $\mathbb{R}^{n}$.

\begin{itemize}
\item[(i)] (Isotropic multipolar potential) Let $V(x)=\sum_{j=1}^{m}%
\frac{\lambda_{j}}{|x-x^{j}|^{2}}$ with $x^{j}=(x_{1}^{j},x_{2}^{j}%
,...,x_{n}^{j})$ and $\lambda_{j}\in\mathbb{R}$. The condition (\ref{pot-cond}%
) is verified for
\begin{equation}
\sum_{j=1}^{m}\left\vert \lambda_{j}\right\vert <(k-2)(n-k).\label{cond-aux-3}%
\end{equation}
The better restriction in (\ref{cond-aux-3}) holds for $k=\frac{n+2}{2}$. In
this case we obtain%
\begin{equation}
\sum_{j=1}^{m}\left\vert \lambda_{j}\right\vert <\lambda_{\ast}=\frac
{(n-2)^{2}}{4}.\label{cond-aux-3-2}%
\end{equation}

\item[(ii)] (Dipole potential) Let $V(x)=\frac{d\cdot x}{|x|^{3}}$ where
$x=(x_{1},x_{2},...,x_{n})$ and $d=(d_{1},d_{2},...,d_{n})$. The condition
(\ref{pot-cond}) is satisfied for
\begin{equation}
\left\vert d\right\vert <\frac{\pi}{(n-2)}\frac{(k-2)(n-k)}{\beta\left(
\frac{1}{2},\frac{n-1}{2}\right)  }.\label{cond-aux-4}%
\end{equation}
For $k=\frac{n+2}{2}$ one obtains
\begin{equation}
|d|<\frac{\pi}{(n-2)}\frac{\lambda_{\ast}}{\beta\left(  \frac{1}{2},\frac
{n-1}{2}\right)  },\label{cond-aux-4-2}%
\end{equation}
which corresponds to the maximum of the R.H.S. in (\ref{cond-aux-4}).

\item[(iii)] (Anisotropic multipolar potential) Let $V(x)=\sum_{j=1}^{m}%
\frac{(x-x^{j}).d^{j}}{|x-x^{j}|^{3}}$ with $x^{j}=(x_{1}^{j},x_{2}%
^{j},...,x_{n}^{j})$ and $d^{j}=(d_{1}^{j},d_{2}^{j},...,d_{n}^{j})$. The
condition (\ref{pot-cond}) is verified for%
\begin{equation}
\sum_{j=1}^{m}\left\vert d^{j}\right\vert <\frac{\pi}{(n-2)}\frac
{(k-2)(n-k)}{\beta\left(  \frac{1}{2},\frac{n-1}{2}\right)  }%
.\label{cond-aux-5}%
\end{equation}

Similarly to item (ii), the maximum of the R.H.S. of (\ref{cond-aux-5}) is
achieved at $k=\frac{n+2}{2}.$ In this case, we get%
\begin{equation}
\sum_{j=1}^{m}\left\vert d^{j}\right\vert <\frac{\pi}{(n-2)}\frac
{\lambda_{\ast}}{\beta\left(  \frac{1}{2},\frac{n-1}{2}\right)  }%
.\label{cond-aux-5-2}%
\end{equation}

\end{itemize}
\end{corollary}

\bigskip

Assuming a certain homogeneity on $u_{0}$ and $V$, then the solution $u$ is
self-similar. The solution is positive when $V$ and $u_{0}$ are also.
Depending on the radial symmetry of $u_{0}$ and $V$, we also investigate if
$u$ is radially symmetric or not.

\begin{theorem}
\label{t.ii} Under the hypotheses of Theorem \ref{t.i}.

\begin{itemize}
\item[(i)] (Self-similarity) Assume that $u_{0}$ and $V$ are homogeneous of
degree $-(n-k)$ and $-2$, respectively. Then the solution $u$ satisfies
\[
u_{\lambda}(x,t)=\lambda^{n-k}u(\lambda x,\lambda^{2}t),
\]
for all $x\in\mathbb{R}^{n}$ and $t>0,$ i.e. $u$ is self-similar.

\item[(ii)] (Positivity) If $V,u_{0}\geq0$ (resp. $V\geq0$, $u_{0}\leq0$) and
$u_{0}\not \equiv 0$, then $u$ is positive (resp. negative).

\item[(iii)] (Radial symmetry) Let $V$ be radially symmetric. The solution $u$
obtained in Theorem \ref{t.i} is radially symmetric for each $t>0$ if and only
if $u_{0}$ is radially symmetric.
\end{itemize}
\end{theorem}

\bigskip

For $V$ as in (\ref{pot1}) and $0\leq\left\vert \lambda\right\vert
<\lambda_{\ast}$, the problem (\ref{heat1})-(\ref{heat3}) (and then
(\ref{mild})) has the explicit stationary solutions
\begin{equation}
\omega_{1}=A_{1}\left\vert x\right\vert ^{-\frac{n-2}{2}+l}\text{ and }%
\omega_{2}=A_{2}\left\vert x\right\vert ^{-\frac{n-2}{2}-l}, \label{stat2}%
\end{equation}
where $A_{i}$'s are arbitrary real constant and $l=\sqrt{\lambda_{\ast
}-\lambda}$ (see \cite{Vazquez1}). A direct computation shows that the
indexes
\begin{equation}
k_{1}=\frac{n+2}{2}+l\text{ and }k_{2}=\frac{n+2}{2}-l \label{ki}%
\end{equation}
are the unique ones such that $\omega_{1}\in PM^{k_{1}}$ and $\omega_{2}\in
PM^{k_{2}}.$ The uniqueness assertion in Theorem \ref{t.i} says that the
unique solution of (\ref{mild}) in $X_{k_{i}}$, with initial data $\omega_{i}%
$, is the corresponding stationary one.

In the next result we analyze the asymptotic behavior of solutions and give in
particular a criterion for their converge towards a stationary state. Here we
employ ideas introduced in \cite{Cannone1}.

\begin{theorem}
\label{t.iii} Assume the hypotheses of Theorem \ref{t.i}. Let $u$ and $v$ be
two solutions for (\ref{mild}) in $X_{k}$ corresponding to the data $u_{0}$
and $v_{0}$ $\in PM^{k}$, respectively. We have that
\begin{equation}
\lim_{t\rightarrow\infty}\left\Vert u(\cdot,t)-v(\cdot,t)\right\Vert _{PM^{k}%
}=0 \label{asymp1}%
\end{equation}
if and only if
\begin{equation}
\lim_{t\rightarrow\infty}\left\Vert G(t)(u_{0}-v_{0})\right\Vert _{PM^{k}}=0.
\label{cond1}%
\end{equation}

The condition (\ref{cond1}) is verified for $u_{0}=v_{0}+\varphi$ with
$\varphi\in\mathcal{S}(\mathbb{R}^{n}).$ In particular, if $V(x)=\frac
{\lambda}{|x|^{2}}$ with $|\lambda|<(k_{i}-2)(n-k_{i}),$ where $k_{i}$ is as
in (\ref{ki}), and $u_{0}=\omega_{i}+\varphi$ with $\varphi\in\mathcal{S}%
(\mathbb{R}^{n})$ then
\begin{equation}
\lim_{t\rightarrow\infty}\left\Vert u(\cdot,t)-\omega_{i}\right\Vert
_{PM^{k_{i}}}=0, \label{asymp2}%
\end{equation}
which shows an attractor-basin around each stationary solution $\omega_{i}$ in
$PM^{k_{i}}$.
\end{theorem}

\bigskip

\begin{remark}
\label{Rem-Asymp}Let $V(x)=\lambda|x|^{-2}$ and $k$ be as in Theorem \ref{t.i}
and $u_{0}$ a non-radial homogeneous function of degree $-(n-k),$ for instance
$u_{0}=$ $x_{j}\left\vert x\right\vert ^{-(n-k+1)}$. It follows from Theorem
\ref{t.ii} (i) that the corresponding solution $u$ is self-similar and a
simple computation shows that it is not stationary. In the case $k=k_{i},$ the
asymptotic behavior of $u$ is not described by the stationary solution
$\omega_{i}$, because $\psi=u_{0}-\omega_{i}\not \equiv 0$ is homogeneous of
degree $-(n-k_{i})$ and then
\[
\lim_{t\rightarrow\infty}\left\Vert G(t)\psi\right\Vert _{PM^{k_{i}}%
}=\left\Vert G(1)\psi\right\Vert _{PM^{k_{i}}}\neq0.
\]
Instead, we obtain from (\ref{asymp1}) a basin of attraction around the
self-similar solution $u$. More precisely, if $v_{0}=u_{0}+\varphi$ then t{he
perturbed solution }$v$ is attracted to the solution $u$ in the sense of
(\ref{asymp1}). Indeed, (\ref{cond1}) induces an equivalent relation in the
set of initial data $PM^{k}$, that is, $u_{0}\sim v_{0}$ if and only if we
have (\ref{cond1}).

The above considerations show a diversified asymptotic behavior of solutions
in $PM^{k}$ with infinitely many possible asymptotics which are characterized
by equivalent classes of initial data.
\end{remark}

The outline of this paper is as follows. In Section 2.1 we recall some basic
properties about Fourier transform and convolution operators useful to handle
(\ref{mild}), and give estimates for the operators (\ref{Group}) and
(\ref{Tv}). Theorems \ref{t.i} and \ref{t.ii} and Corollary \ref{c.i} are
proved in subsections 2.2, 2.3 and 2.4. Subsection 2.5 is devoted to proving
Theorem \ref{t.iii}.

\section{Proof of the Results}

\subsection{\bigskip Basic estimates in $PM^{a}$}

We start by recalling some results about convolution and Fourier transform of
homogeneous functions which will be useful to perform estimates with explicit
constants in $PM^{k}$-spaces (see \cite[p. 124]{Li-Lo} and \cite[p.
160]{Stein}, respectively).

\begin{lemma}
\label{conv} Let $0<\theta_{1}<n$, $0<\theta_{2}<n$ and $0<\theta_{1}%
+\theta_{2}<n$. Then
\begin{equation}
(|x|^{\theta_{1}-n}\ast|x|^{\theta_{2}-n})(y)=\int_{\mathbb{R}^{n}}%
|z|^{\theta_{1}-n}|y-z|^{\theta_{2}-n}dz=K(\theta_{1},\theta_{2}%
,n)|y|^{\theta_{1}+\theta_{2}-n}, \label{e.integral}%
\end{equation}
where $K(\theta_{1},\theta_{2},n)=(\nu_{\theta_{1}}\nu_{\theta_{2}}%
\nu_{n-\theta_{1}-\theta_{2}})/(\nu_{\theta_{1}+\theta_{2}}\nu_{n-\theta_{1}%
}\nu_{n-\theta_{2}})$ and $\nu_{\theta}=\pi^{-\theta/2}\Gamma(\theta/2)$.
\end{lemma}

\begin{lemma}
\label{fourier} Suppose that $\alpha$ is a complex number such that
$0<Re(\alpha)<n$ and $P_{l}(x)$ is a harmonic polynomial on $\mathbb{R}^{n}$
homogeneous of degree $l$. If $K(x)=\frac{P_{l}(x)}{|x|^{n+l-\alpha}}$ then
$\widehat{K}(\xi)=\gamma_{l,\alpha}\frac{P_{l}(\xi)}{|\xi|^{l+\alpha}}$,
where
\begin{equation}
\gamma_{l,\alpha}=\frac{i^{-l}\pi^{(n/2)-\alpha}\Gamma\left(  \frac{l+\alpha
}{2}\right)  }{\Gamma\left(  \frac{n+l-\alpha}{2}\right)  }. \label{const-1}%
\end{equation}

\end{lemma}

In the following we recall an estimate in $PM^{k}$-spaces for the heat
semigroup (\ref{Group}) (see e.g. \cite{Cannone1}). The proof is included for
the reader convenience.

\begin{lemma}
\label{boundB1} If $f\in PM^{k}$ with $k\geq0$ then $G(t)f\in BC_{w}\left(
[0,\infty);PM^{k}\right)  $ and%
\begin{equation}
\sup_{t>0}\Vert G(t)f\Vert_{PM^{k}}\leq\Vert f\Vert_{PM^{k}}\text{.}\label{l1}%
\end{equation}

\end{lemma}

\textbf{Proof.} \ We have that

\begin{equation}
\left\vert \xi\right\vert ^{k}\left\vert \widehat{G(t)f}(\xi)\right\vert
\leq\left\vert \xi\right\vert ^{k}|e^{-4\pi^{2}t|\xi|^{2}}\hat{f}(\xi
)|\leq\left\vert \xi\right\vert ^{k}|\hat{f}(\xi)|, \label{l1-1}%
\end{equation}
which gives (\ref{l1}) after applying $ess\sup_{x\in\mathbb{R}^{n}}$ in both
sides of (\ref{l1-1}). \qed

\bigskip

The next lemma gives an estimate for the linear operator (\ref{Tv}) in
$PM^{k}$-spaces by expliciting the dependence of the norm of $V$ and giving an
exact expression for the constant.

\begin{lemma}
\label{boundB2} Let $0<b_{1},b_{2}<n$ be such that $n<b_{1}+b_{2}<2n$. For
$K(\theta_{1},\theta_{2},n)$ as in Lemma \ref{e.integral}, we set
\begin{equation}
C_{b_{1},b_{2}}=\frac{1}{4\pi^{2}}K(n-b_{1},n-b_{2},n).\label{const-aux-1}%
\end{equation}
Then $L_{V}(u)\in BC_{w}\left(  [0,\infty);PM^{b}\right)  $ with
$b=b_{1}+b_{2}+2-n$ and
\begin{equation}
\sup_{t>0}\Vert L_{V}(u)(t)\Vert_{PM^{b}}\leq C_{b_{1},b_{2}}\Vert
V\Vert_{PM^{b_{1}}}\sup_{t>0}\Vert u(\cdot,t)\Vert_{PM^{b_{2}}},\label{l2}%
\end{equation}
for all $V\in PM^{b_{1}}$ and $u\in X_{b_{2}}.$
\end{lemma}

\textbf{Proof.}{} Using Lemma \ref{e.integral} we obtain
\begin{align}
|\widehat{V}\ast\widehat{u}(\xi)|  &  \leq\int_{\mathbb{R}^{n}}\left\vert
\widehat{V}(\xi-\eta)\widehat{u}(\eta)\right\vert d\eta\nonumber\\
&  \leq\int_{\mathbb{R}^{n}}\frac{1}{|\xi-\eta|^{b_{1}}}\frac{1}{|\eta
|^{b_{2}}}d\eta\left\Vert V\right\Vert _{PM^{b_{1}}}\left\Vert u\right\Vert
_{PM^{b_{2}}}\nonumber\\
&  \leq K(n-b_{1},n-b_{2},n)\frac{1}{|\xi|^{b_{1}+b_{2}-n}}\left\Vert
V\right\Vert _{PM^{b_{1}}}\left\Vert u\right\Vert _{PM^{b_{2}}}. \label{l4}%
\end{align}
It follows from (\ref{Tv}) and (\ref{l4}) that
\begin{align*}
\left\vert \widehat{L_{V}(u)}(\xi)\right\vert  &  \leq\int_{0}^{t}e^{-4\pi
^{2}(t-s)|\xi|^{2}}|\widehat{V}\ast\left(  \widehat{u}(\xi,s)\right)  |ds\\
&  \leq\int_{0}^{t}e^{-4\pi^{2}(t-s)|\xi|^{2}}\frac{K(n-b_{1},n-b_{2},n)}%
{|\xi|^{b_{1}+b_{2}-n}}\left\Vert V\right\Vert _{PM^{b_{1}}}\left\Vert
u(\cdot,s)\right\Vert _{PM^{b_{2}}}ds\\
&  \leq\frac{K(n-b_{1},n-b_{2},n)}{|\xi|^{b_{1}+b_{2}-n}}\int_{0}^{t}%
e^{-4\pi^{2}(t-s)|\xi|^{2}}ds\left\Vert V\right\Vert _{PM^{b_{1}}}\sup
_{t>0}\left\Vert u(\cdot,t)\right\Vert _{PM^{b_{2}}}\\
&  \leq\frac{1}{4\pi^{2}}\frac{K(n-b_{1},n-b_{2},n)}{|\xi|^{b_{1}+b_{2}+2-n}%
}(1-e^{-4\pi^{2}t|\xi|^{2}})\left\Vert V\right\Vert _{PM^{b_{1}}}\sup
_{t>0}\left\Vert u(\cdot,t)\right\Vert _{PM^{b_{2}}},
\end{align*}
which yields the desired inequality.

The weak-time continuity in $\mathcal{S^{\prime}}(\mathbb{R}^{n})$ is left to
the reader (see e.g. \cite{Cannone1} and \cite{Yamazaki}). \qed

\subsection{Proof of Theorem \ref{t.i}}

\textbf{Part (i): }Let us set
\[
\tau=C_{n-2,k}\left\Vert V\right\Vert _{PM^{n-2}}.
\]
Lemma \ref{boundB2} with $(b_{1},b_{2})=(n-2,k)$ yields%

\begin{align}
\Vert L_{V}(u)(t)-L_{V}(v)(t)\Vert_{X_{k}} &  =\sup_{t>0}\Vert L_{V}%
(u-v)(t)\Vert_{PM^{k}}\nonumber\\
&  \leq C_{n-2,k}\left\Vert V\right\Vert _{PM^{n-2}}\Vert u-v\Vert_{X_{k}%
},\label{b11}%
\end{align}
and then
\[
\left\Vert L_{V}\right\Vert _{X_{k}\rightarrow X_{k}}\leq\tau<1.
\]
Also, Lemma \ref{boundB1} gives us
\begin{equation}
\Vert G(t)u_{0}\Vert_{X_{k}}=\sup_{t>0}\Vert G(t)u_{0}\Vert_{PM^{k}}\leq\Vert
u_{0}\Vert_{PM^{k}}.\label{b22}%
\end{equation}
The estimate (\ref{b11}) with $v=0$ and (\ref{b22}) imply that the operator
$H:X_{k}\rightarrow X_{k}$ such that $H(u)=G(t)u_{0}+L_{V}(u)(t)$ is well
defined. Also, we have the estimate
\begin{equation}
\Vert H(u)-H(v)\Vert_{X_{k}}\leq\Vert L_{V}(u)(t)-L_{V}(v)(t)\Vert_{X_{k}}%
\leq\tau\Vert u-v\Vert_{X_{k}},\text{ for all }u,v\in X_{k},\label{ineq-aux-1}%
\end{equation}
which shows that $H$ is a contraction in $X_{k}$. Now the Banach fixed point
theorem assures the existence of a unique solution $u\in X_{k}$ for
(\ref{mild}). \qed

\textbf{Part (ii): } From (\ref{const-aux-1}) and Lemma \ref{conv}, we can
compute $C_{n-2,k}$ explicitly as
\begin{align}
C_{n-2,k} &  =\frac{K(2,n-k,n)}{4\pi^{2}}\nonumber\\
&  =\frac{\pi^{n/2}\Gamma(1)\Gamma\left(  \frac{(n-k)}{2}\right)
\Gamma\left(  \frac{(k-2)}{2}\right)  }{4\pi^{2}\Gamma\left(  \frac{(n-2)}%
{2}\right)  \Gamma\left(  \frac{k}{2}\right)  \Gamma\left(  \frac{(2+n-k)}%
{2}\right)  }\nonumber\\
&  =\frac{\pi^{n/2}(n-2)}{2\pi^{2}\Gamma\left(  \frac{n}{2}\right)
(k-2)(n-k)}.\label{est.c.i.1}%
\end{align}
It follows from Theorem \ref{t.i} (i) that
\begin{equation}
\Vert V\Vert_{PM^{n-2}}<\frac{1}{C_{n-2,k}}=\frac{2\pi^{2}\Gamma\left(
\frac{n}{2}\right)  (k-2)(n-k)}{\pi^{n/2}(n-2)}.\label{est.V}%
\end{equation}
For $V=\frac{\lambda}{|x|^{2}}$, applying Lemma \ref{fourier} with
$\alpha=n-2$ and $l=0$, we obtain
\[
\widehat{V}(\xi)=\lambda\pi^{2-\frac{n}{2}}\Gamma\left(  \frac{n-2}{2}\right)
|\xi|^{2-n},
\]
which gives $V\in PM^{n-2}$ with
\begin{equation}
\Vert V\Vert_{PM^{n-2}}=\left\vert \lambda\right\vert \pi^{2-\frac{n}{2}%
}\Gamma\left(  \frac{n-2}{2}\right)  .\label{norm-V1}%
\end{equation}
In view of (\ref{norm-V1}), the condition (\ref{est.V}) can be expressed by
means of the size of $\left\vert \lambda\right\vert $ as
\begin{align*}
|\lambda| &  <\frac{1}{\pi^{2-\frac{n}{2}}\Gamma\left(  \frac{n-2}{2}\right)
C_{n-2,k}}\\
&  =\frac{2\pi^{2}\Gamma\left(  \frac{n}{2}\right)  (n-k)(k-2)}{\pi^{2}%
\Gamma\left(  \frac{n-2}{2}\right)  (n-2)}=(n-k)(k-2).
\end{align*}
\qed

\textbf{Part (iii): }Let $u$ and $v$ be two solutions obtained from
\textit{item (i)} with data $V,u_{0}$ and $W,v_{0}$, respectively. Firstly,
using (\ref{ineq-aux-1}), we obtain
\[
\left\Vert v\right\Vert _{X_{k}}=\left\Vert H(v)\right\Vert _{X_{k}}%
\leq\left\Vert G(t)v_{0}\right\Vert _{X_{k}}+\left\Vert L_{W}(v)\right\Vert
_{X_{k}}\leq\Vert v_{0}\Vert_{PM^{k}}+C_{n-2,k}\left\Vert W\right\Vert
_{PM^{n-2}}\left\Vert v\right\Vert _{X_{k}},
\]
which implies
\begin{equation}
\Vert v\Vert_{X_{k}}\leq\frac{\Vert v_{0}\Vert_{PM^{k}}}{1-C_{n-2,k}\left\Vert
W\right\Vert _{PM^{n-2}}}.\text{ } \label{aux-proof-1}%
\end{equation}
Subtracting the respective equations verified by $u,v$, and afterwards
applying the norm $\Vert\cdot\Vert_{X_{k}}$, we estimate
\begin{align}
\left\Vert u-v\right\Vert _{X_{k}}  &  =\left\Vert G(t)(u_{0}-v_{0}%
)+L_{V}(u-v)+L_{V-W}(v)\right\Vert _{X_{k}}\nonumber\\
&  \leq\Vert u_{0}-v_{0}\Vert_{PM^{k}}+C_{n-2,k}\left\Vert V\right\Vert
_{PM^{n-2}}\left\Vert u-v\right\Vert _{X_{k}}\nonumber\\
&  +C_{n-2,k}\left\Vert V-W\right\Vert _{PM^{n-2}}\left\Vert v\right\Vert
_{X_{k}}. \label{aux-proof-2}%
\end{align}
The estimates (\ref{aux-proof-1}) and (\ref{aux-proof-2}) lead us to
\[
\left\Vert u-v\right\Vert _{X_{k}}\leq\frac{1}{1-C_{n-2,k}\left\Vert
V\right\Vert _{PM^{n-2}}}\left(  \Vert u_{0}-v_{0}\Vert_{PM^{k}}%
+\frac{C_{n-2,k}\Vert v_{0}\Vert_{PM^{k}}}{1-C_{n-2,k}\left\Vert W\right\Vert
_{PM^{n-2}}}\left\Vert V-W\right\Vert _{PM^{n-2}}\right)  ,
\]
as desired. \qed

\subsection{Proof of Corollary \ref{c.i}}

\textbf{Part (i): }(Isotropic multipolar potential) Using the translation
property of Fourier transform and Lemma \ref{fourier}, we obtain
\[
\widehat{V}(\xi)=\sum_{j=1}^{m}\lambda_{j}\pi^{2-\frac{n}{2}}\Gamma\left(
\frac{n-2}{2}\right)  e^{-2\pi i(x^{j}.\xi)}|\xi|^{2-n},
\]
which gives $V\in PM^{n-2}$ and
\begin{align*}
\Vert V\Vert_{PM^{n-2}}  &  \leq\pi^{2-\frac{n}{2}}\Gamma\left(  \frac{n-2}%
{2}\right)  \left(  \Sigma_{j=1}^{m}\left\vert \lambda_{j}\right\vert \right)
\\
&  <\frac{1}{C_{n-2,k}},
\end{align*}
provided that $\Sigma_{j=1}^{m}\left\vert \lambda_{j}\right\vert <(n-k)(k-2)$.

\bigskip

\textbf{Part (ii):} (Dipole potential) Lemma \ref{fourier} with $\alpha=n-2$
and $l=1$ yields
\begin{align*}
\widehat{V}(\xi)  &  =\left(  \sum_{j=1}^{m}\frac{d_{j}x_{j}}{|x|^{3}}\right)
^{\wedge}=\sum_{j=1}^{m}2i^{-1}\pi^{\frac{3-n}{2}}\Gamma\left(  \frac{n-1}%
{2}\right)  \frac{d_{j}\xi_{j}}{|\xi|^{n-1}}\\
&  =2i^{-1}\pi^{\frac{3-n}{2}}\Gamma\left(  \frac{n-1}{2}\right)  \frac
{d\cdot\xi}{|\xi|^{n-1}}.
\end{align*}
It follows that $V\in PM^{n-2}$ and
\[
\Vert V\Vert_{PM^{n-2}}\leq2\pi^{\frac{3-n}{2}}\Gamma\left(  \frac{n-1}%
{2}\right)  \left\vert d\right\vert .
\]
Then, the condition (\ref{est.V}) (i.e. (\ref{pot-cond})) is verified when
\begin{align}
|d|  &  <\frac{1}{2\pi^{\frac{3-n}{2}}\Gamma\left(  \frac{n-1}{2}\right)
}\frac{2\pi^{2-n/2}(n-k)(k-2)\Gamma\left(  \frac{n}{2}\right)  }%
{(n-2)}\nonumber\\
&  =\frac{\pi^{1/2}(n-k)(k-2)\Gamma\left(  \frac{n}{2}\right)  }%
{(n-2)\Gamma\left(  \frac{n-1}{2}\right)  }\nonumber\\
&  =\frac{\pi(n-k)(k-2)}{(n-2)}\frac{1}{\beta\left(  \frac{1}{2},\frac{n-1}%
{2}\right)  }. \label{aux-proof-4}%
\end{align}

\textbf{Part (iii): }(Anisotropic multipolar potential) Computing the Fourier
transform, it follows that
\[
\widehat{V}(\xi)=\sum_{j=1}^{m}2i^{-1}\pi^{\frac{3-n}{2}}\Gamma\left(
\frac{n-1}{2}\right)  e^{-2\pi i(x^{j}.\xi)}\frac{\xi.d^{j}}{|\xi|^{n-1}},
\]
and then, similarly to\textit{ item (ii)}, the condition (\ref{pot-cond}) for
$V$ in Theorem \ref{t.i} is verified for
\[
\sum_{j=1}^{m}\left\vert d^{j}\right\vert <\frac{\pi(n-k)(k-2)}{(n-2)}\frac
{1}{\beta\left(  \frac{1}{2},\frac{n-1}{2}\right)  }.
\]
\qed

\subsection{Proof of Theorem \ref{t.ii}}

\textbf{Part (i): }Let $u$ be the solution given in Theorem \ref{t.i} (i). It
is not difficult to check that $u_{\lambda}(x,t)=\lambda^{n-k}u(\lambda
x,\lambda^{2}t)\in X_{k}$ also verifies (\ref{mild}) when $u_{0}(x)$ and
$V(x)$ are homogeneous of degree $-(n-k)$ and $-2$, respectively. Now the
uniqueness statement in Theorem \ref{t.i} gives us $u=u_{\lambda}$ for all
$\lambda>0$, as required.\noindent

\bigskip

\textbf{Part (ii): }We should prove only the positivity statement, because the
proof of the one between brackets is similar. Recall first that $F\in
\mathcal{S}^{\prime}(\mathbb{R}^{n})$ is said to be nonnegative (resp.
nonpositive) if $\left\langle F,\varphi\right\rangle \geq0$ (resp. $\leq0$),
for all $\varphi\geq0$ and $\varphi\in\mathcal{S}(\mathbb{R}^{n})$. Also, $F$
is positive (resp. negative) when $\left\langle F,\varphi\right\rangle >0$
(resp. $<0$), for all $\varphi>0$ and $\varphi\in\mathcal{S}(\mathbb{R}^{n})$.

Note that $u_{1}=G(t)u_{0}$ is a positive distribution in $\mathcal{S}%
^{\prime}(\mathbb{R}^{n})$, for $t>0$, when $u_{0}$ is nonnegative and
$u_{0}\not \equiv 0$. Since the solution $u$ has been obtained via Banach
fixed point theorem, it is the limit of the Picard interaction
\begin{equation}
u_{1}=G(t)u_{0}\text{ \ and }u_{b+1}=u_{1}+L_{V}(u_{b})\text{, \ }%
b\in\mathbb{N}.\label{sequencee}%
\end{equation}
An induction argument shows that all elements of (\ref{sequencee}) are
positive distribution in $\mathcal{S}^{\prime}(\mathbb{R}^{n})$, for $t>0$.
Since $u_{b}\rightarrow u$ in $X_{k}$, we have that $u_{b}\rightarrow u$ in
$\mathcal{S}^{\prime}(\mathbb{R}^{n}),$ for $t>0.$ It follows that
$u(\cdot,t)$ is a nonnegative distribution, for $t>0,$ because the convergence
in $\mathcal{S}^{\prime}(\mathbb{R}^{n})$ preserves nonnegativity. As $u_{1}$
is positive and $L_{V}(u)$ is nonnegative, it follows that
\[
\left\langle u(\cdot,t),\varphi\right\rangle =\left\langle u_{1}%
(\cdot,t),\varphi\right\rangle +\left\langle L_{V}(u)(t),\varphi\right\rangle
\geq\left\langle u_{1}(\cdot,t),\varphi\right\rangle >0,\text{ for }t>0\text{,
}%
\]
for all $\varphi>0$ and $\varphi\in\mathcal{S}(\mathbb{R}^{n}).$

\noindent

\textbf{Part (iii): }Let $u_{0}$ and $V$ be radially symmetric. As the heat
flow preserves radial symmetry, \ it follows that $u_{1}=G(t)u_{0}$ is
radially symmetric, for each fixed $t>0$. Also, $L_{V}(u)$ is radially
symmetric provided that $u$ is also radially symmetric. One can prove by
induction that $\{u_{b}\}_{b\geq1}$ (see (\ref{sequencee})) is radially
symmetric, for each fixed $t>0.$ Since $u_{b}\rightarrow u$ in $X_{k}$ and
Fourier transform preserves radial symmetry, we get that $u$ is radially
symmetric, for each fixed $t>0.$

Assume now that $u_{0}$ is not radially symmetric and $V$ is radially
symmetric. Suppose, to the contrary, that $u$ were radially symmetric, then
$L_{V}(u)$ also would be radially symmetric. So, $G(t)u_{0}$ $=u-L_{V}(u)$
would be radially symmetric, which gives a contradiction because
$(G(t)u_{0})^{\wedge}=e^{-\left\vert \xi\right\vert ^{2}t}\widehat{u}_{0}$ is
radially symmetric if and only if $\widehat{u}_{0}$ is radially symmetric. \qed

\subsection{Proof of Theorem \ref{t.iii}}

We only prove that (\ref{cond1}) implies (\ref{asymp1}). The converse
statement follows similarly and is left to the reader. Subtracting the
equations satisfied by $u$ and $v$, and afterwards computing the $PM^{k}%
$-norm, we obtain
\begin{equation}
\left\Vert u(\cdot,t)-v(\cdot,t)\right\Vert _{PM^{k}}\leq\left\Vert
G(t)(u_{0}-v_{0})\right\Vert _{PM^{k}}+J_{1}(t)+J_{2}(t) \label{ineq11}%
\end{equation}
where
\begin{align*}
J_{1}(t)  &  =4\pi^{2}C_{n-2,k}\left\Vert V\right\Vert _{PM^{n-2}}\sup_{\xi
\in\mathbb{R}^{n}}\int_{0}^{\delta t}|\xi|^{2}e^{-4\pi^{2}\left\vert
\xi\right\vert ^{2}(t-s)}\Vert u(\cdot,s)-v(\cdot,s)\Vert_{PM^{k}}ds\\
J_{2}(t)  &  =4\pi^{2}C_{n-2,k}\left\Vert V\right\Vert _{PM^{n-2}}\sup_{\xi
\in\mathbb{R}^{n}}\int_{\delta t}^{t}|\xi|^{2}e^{-4\pi^{2}\left\vert
\xi\right\vert ^{2}(t-s)}\Vert u(\cdot,s)-v(\cdot,s)\Vert_{PM^{k}}ds.
\end{align*}
with $\delta>0$ being a constant that will be chosen later. Using that
\[
\sup_{\xi\in\mathbb{R}^{n}}t|\xi|^{2}e^{-t(1-s)4\pi^{2}|\xi|^{2}}=\frac
{e^{-1}}{4\pi^{2}(1-s)},
\]
and the change $s=tz$ in $J_{1}(t)$, we estimate
\begin{align}
J_{1}(t)  &  \leq4\pi^{2}C_{n-2,k}\left\Vert V\right\Vert _{PM^{n-2}}\sup
_{\xi\in\mathbb{R}^{n}}\int_{0}^{\delta}t|\xi|^{2}e^{-t(1-s)4\pi^{2}|\xi|^{2}%
}\Vert u(\cdot,ts)-v(\cdot,ts)\Vert_{PM^{k}}ds\nonumber\\
&  \leq C\int_{0}^{\delta}(1-s)^{-1}\Vert u(\cdot,ts)-v(\cdot,ts)\Vert
_{PM^{k}}ds. \label{A1}%
\end{align}
The term $J_{2}(t)$ can be estimated directly by
\begin{align}
J_{2}(t)  &  \leq4\pi^{2}C_{n-2,k}\left\Vert V\right\Vert _{PM^{n-2}}\left(
\sup_{\xi\in\mathbb{R}^{n}}\int_{\delta t}^{t}|\xi|^{2}e^{-(t-s)4\pi^{2}%
|\xi|^{2}}ds\right)  \left(  \sup_{\delta t<s<t}\Vert u(\cdot,s)-v(\cdot
,s)\Vert_{PM^{k}}\right) \nonumber\\
&  =C_{n-2,k}\left\Vert V\right\Vert _{PM^{n-2}}\sup_{\delta t<s<t}\Vert
u(\cdot,s)-v(\cdot,s)\Vert_{PM^{k}}, \label{A2}%
\end{align}
because $\int_{\delta t}^{t}|\xi|^{2}e^{-(t-s)4\pi^{2}|\xi|^{2}}ds=\frac
{1}{4\pi^{2}}\left(  1-e^{-4\pi^{2}(1-\delta)t|\xi|^{2}}\right)  $. Noting
that
\[
\Gamma=\limsup_{t\rightarrow\infty}\Vert u(\cdot,t)-v(\cdot,t)\Vert_{PM^{k}%
}\leq(\left\Vert u\right\Vert _{X_{k}}+\left\Vert v\right\Vert _{X_{k}%
})<\infty,
\]
we can calculate the superior limit in (\ref{ineq11}), and then use (\ref{A1})
and (\ref{A2}) in order to obtain
\[
\Gamma\leq\left(  C\log\left(  \frac{1}{1-\delta}\right)  +C_{n-2,k}\left\Vert
V\right\Vert _{PM^{n-2}}\right)  \Gamma=M\Gamma.
\]
In view of $C_{n-2,k}\left\Vert V\right\Vert _{PM^{n-2}}<1$ (see
(\ref{pot-cond})), one can take $\delta>0$ in such a way that $0<M<1$, and so
$\Gamma=0$, as required.

The further conclusions in the statement follow by employing (\ref{asymp1})
with $v(x,t)\equiv\omega_{i}$ and noting that $\lim_{t\rightarrow\infty
}\left\Vert G(t)\varphi\right\Vert _{PM^{k}}=0$ when $\varphi\in
\mathcal{S}\mathbb{(R}^{n}\mathbb{)}.$ \qed


\begin{thebibliography}{99}                                                                                               
\bibitem {Abdellaoui-Peral1}B. Abdellaoui, I. Peral, V. Felli, Existence and
multiplicity for perturbations of an equation involving a Hardy inequality and
the critical Sobolev exponent in the whole of $R^{n}$, Adv. Differential
Equations 9 (2004), 481--508.

\bibitem {Abdellaoui-Peral2}B. Abdellaoui, I. Peral, A. Primo, Strong
regularizing effect of a gradient term in the heat equation with the Hardy
potential, J. Funct. Anal. 258 (4) (2010), 1247--1272.

\bibitem {Abdellaoui-Peral3}B. Abdellaoui, I. Peral, A. Primo, Optimal results
for parabolic problems arising in some physical models with critical growth in
the gradient respect to a Hardy potential, Adv. Math. 225 (6) (2010), 2967--3021.

\bibitem {Abdellaoui-Peral4}B. Abdellaoui, I. Peral, A. Primo, Influence of
the Hardy potential in a semilinear heat equation, Proc. Roy. Soc. Edinburgh
Sect. A 139 (5) (2009), 897--926.

\bibitem {Baras1}P. Baras and J. Goldstein, The heat equation with a singular
potential, Trans. Amer. Math. Soc. 284 (1984), 121--139.

\bibitem {Biler-Karch}P. Biler, M. Cannone, I. A. Guerra, G. Karch, Global
regular and singular solutions for a model of gravitating particles, Math.
Ann. 330 (4) (2004), 693--708.

\bibitem {Cabre-Martel1}X. Cabr\'{e}, Y. Martel, Existence versus explosion
instantan\'{e}e pour des \'{e}quations de la chaleur lin\'{e}aires avec
potentiel singulier, C. R. Acad. Sci. Paris S\'{e}r. I Math. 329 (11) (1999), 973--978.

\bibitem {Cannone1}M. Cannone, G. Karch, Smooth or singular solutions to the
Navier-Stokes system, J. Differential Equations 197 (2004), 247--274.

\bibitem {Car-Fer1}J. A. Carrillo, L.C.F. Ferreira, Self-similar solutions and
large time asymptotics for the dissipative quasi-geostrophic equation,
Monatsh. Math. 151 (2) (2007), 111--142.

\bibitem {Chaudhuri}N. Chaudhuri, K. Sandeep, On a heat problem involving the
perturbed Hardy-Sobolev operator, Proc. Roy. Soc. Edinburgh Sect. A 134 (4)
(2004), 683--693.

\bibitem {Chaves}M. Chaves, J. Garc\'{\i}a Azorero, On bifurcation and
uniqueness results for some semilinear elliptic equations involving a singular
potential, J. Eur. Math. Soc. 8 (2) (2006), 229--242.

\bibitem {Peral1}A. Dall'Aglio, D. Giachetti, I. Peral, Results on parabolic
equations related to some Caffarelli-Kohn-Nirenberg inequalities, SIAM J.
Math. Anal. 36 (3) (2004/05), 691--716.

\bibitem {Davila}J. D\'{a}vila, L. Dupaigne, Comparison results for PDEs with
a singular potential, Proc. Roy. Soc. Edinburgh Sect. A 133 (1) (2003), 61--83.

\bibitem {Felli1}V. Felli, A. Pistoia, Existence of blowing-up solutions for a
nonlinear elliptic equation with Hardy potential and critical growth, Comm.
Partial Differential Equations 31 (2006), 21--56.

\bibitem {FMT1}V. Felli, E. M. Marchini, S. Terracini, On Schr\"{o}dinger
operators with multipolar inverse-square potentials, Journal of Func. Anal.
250 (2007), 265-316.

\bibitem {FMT2}V. Felli, E. M. Marchini, S. Terracini, On Schr\"{o}dinger
operators with multisingular inverse-square anisotropic potentials, Indiana
Univ. Math. J. 58 (2009), 617-676.

\bibitem {Fer-Mesq}L.C.F. Ferreira, C.A.A.S. Mesquita, Existence and
symmetries for elliptic equations with multipolar potentials and polyharmonic
operators, to appear in Indiana University Mathematics Journal (2013).

\bibitem {Frank}W.M. Frank, D.J. Land and R.M. Spector, Singular potentials,
Rev. Modern Physics, 43 (1971), 36--98.

\bibitem {Galaktionov1}V.A. Galaktionov, I. V. Kamotski, On nonexistence of
Baras-Goldstein type for higher-order parabolic equations with singular
potentials, Trans. Amer. Math. Soc. 362 (8) (2010), 4117--4136.

\bibitem {Grafakos}L. Grafakos, Classical and modern Fourier analysis, Pearson
Education, Upper Saddle River, NJ, 2004.

\bibitem {Gkikas1}K.T. Gkikas, Existence and nonexistence of energy solutions
for linear elliptic equations involving Hardy-type potentials, Indiana Univ.
Math. J. 58 (5) (2009), 2317--2345.

\bibitem {Goldstein1}G. R. Goldstein, J. A. Goldstein, A. Rhandi, Kolmogorov
equations perturbed by an inverse-square potential, Discrete Contin. Dyn.
Syst. Ser. S 4 (3) (2011), 623--630.

\bibitem {Goldstein2}J. A. Goldstein, Q. S. Zhang, Linear parabolic equations
with strong singular potentials, Trans. Amer. Math. Soc. 355 (1) (2003), 197--211.

\bibitem {Karachalios}N. I. Karachalios, N. B. Zographopoulos, The semiflow of
a reaction diffusion equation with a singular potential, Manuscripta Math. 130
(1) (2009), 63--91.

\bibitem {Kombe}I. Kombe, The linear heat equation with highly oscillating
potential, Proc. Amer. Math. Soc. 132 (9) (2004), 2683--2691.

\bibitem {Landau}L. D. Landau, E. M. Lifshitz, Quantum Mechanics, Pergamon
Press Ltd., London, 1965.

\bibitem {Le-Jan}Y. Le Jan, A. S. Sznitman, Stochastic cascades and
3-dimensional Navier--Stokes equations, Probab. Theory Related Fields 109
(1997), 343--366.

\bibitem {Levy}J. M. L\'{e}vy-Leblond, Electron capture by polar molecules,
Phys. Rev. 153 (1967), 1--4.

\bibitem {Li-Lo}E. Lieb, M. Loss, Analysis, American Mathematical Society,
Providence, RI, 2001.

\bibitem {Liskevich1}V. Liskevich, A. Shishkov, Z. Sobol, Singular solutions
to the heat equations with nonlinear absorption and Hardy potentials, Commun.
Contemp. Math. 14 (2) (2012), 1250013, 28 pp.

\bibitem {Miao1}C. Miao, B. Yuan, Solutions to some nonlinear parabolic
equations in pseudomeasure spaces, Math. Nachr. 280 (1-2) (2007), 171--186.

\bibitem {Peral-Vazquez}I. Peral, J. L. V\'{a}zquez, On the stability or
instability of the singular solution of the semilinear heat equation with
exponential reaction term, Arch. Rational Mech. Anal. 129 (3) (1995), 201--224.

\bibitem {Reyes1}G. Reyes, A. Tesei, Self-similar solutions of a semilinear
parabolic equation with inverse-square potential, J. Differential Equations
219 (1) (2005), 40--77.

\bibitem {Smets}D. Smets, Nonlinear Schr\"{o}dinger equations with Hardy
potential and critical nonlinearities, Trans. Amer. Math. Soc. 357 (2005), 2909--2938

\bibitem {Stein}E. M. Stein, G. Weiss, Introduction to Fourier analysis on
Euclidean spaces, Princeton Mathematical Series 32, Princeton University
Press, Princeton, N.J., 1971.

\bibitem {Vancostenoble1}J. Vancostenoble, Lipschitz stability in inverse
source problems for singular parabolic equations, Comm. Partial Differential
Equations 36 (8) (2011), 1287--1317.

\bibitem {Vazquez1}J. L. Vazquez, E. Zuazua, The Hardy inequality and the
asymptotic behaviour of the heat equation with an inverse-square potential, J.
Funct. Anal. 173 (1) (2000), 103--153.

\bibitem {Yamazaki}M. Yamazaki, The Navier-Stokes equations in the
weak-$L^{n}$ space with time-dependent external force, Math. Ann. 317 (4)
(2000), 635-675.
\end{thebibliography}
\end{document}